\documentclass[12pt]{article}%
\usepackage{graphicx}
\usepackage[intlimits]{amsmath}
\usepackage{latexsym}
\usepackage{amsfonts}
\usepackage{amssymb}%
\makeatletter
\newfont{\footsc}{cmcsc10 at 8truept}
\newfont{\footbf}{cmbx10 at 8truept}
\newfont{\footrm}{cmr10 at 10truept}
\newtheorem{theorem}{Theorem}

\newtheorem{proposition}[theorem]{Proposition}

\newenvironment{proof}[1][Proof]{\noindent{\textbf {#1}  }}  {\hfill$\Box$\bigskip}

\begin{document}

\title{Revisiting two classical results on graph spectra}
\author{Vladimir Nikiforov\\Department of Mathematical Sciences, University of Memphis, \\Memphis TN 38152, USA}
\maketitle

\begin{abstract}
Let $\mu\left(  G\right)  $ and $\mu_{\min}\left(  G\right)  $ be the largest
and smallest eigenvalues of the adjacency matrix of a graph $G$. Our main
results are:

(i) If $H$ is a proper subgraph of a connected graph $G$ of order $n$ and
diameter $D,$ then%
\[
\mu\left(  G\right)  -\mu\left(  H\right)  >\frac{1}{\mu^{2D}\left(  G\right)
n}.
\]

(ii) If $G$ is a connected nonbipartite graph of order $n$ and diameter $D$,
then
\[
\mu\left(  G\right)  +\mu_{\min}\left(  G\right)  >\frac{2}{\mu^{2D}\left(
G\right)  n}.
\]

These bounds have the correct order of magnitude for large $\mu$ and $D$.


\end{abstract}

\section{Introduction}

Our notation is standard (e.g., see \cite{Bol98}, \cite{CDS80}, and
\cite{HoJo88}). In particular, unless specified otherwise, all graphs are
defined on the vertex set $\left[  n\right]  =\left\{  1,...,n\right\}  $ and
$\mu\left(  G\right)  $ and $\mu_{\min}\left(  G\right)  $ stand for the
largest and smallest eigenvalues of the adjacency matrix of a graph $G$.

The aim of this note is to refine quantitatively two well-known results on
graph spectra. The first one, following from Frobenius's theorem on
nonnegative matrices, asserts that if $H$ is a proper subgraph of a connected
graph $G,$ then $\mu\left(  G\right)  >\mu\left(  H\right)  .$ The second one,
due to H. Sachs \cite{Sac64}, asserts that if $G$ is a connected nonbipartite
graph, then $\mu\left(  G\right)  >-\mu_{\min}\left(  G\right)  .$

Our main result is the following theorem.

\begin{theorem}
\label{th1}If $H$ is a proper subgraph of a connected graph $G$ of order $n$
and diameter $D,$ then%
\begin{equation}
\mu\left(  G\right)  -\mu\left(  H\right)  >\frac{1}{\mu^{2D}\left(  G\right)
n}. \label{in1}%
\end{equation}

\end{theorem}

It can be shown that, for large $\mu$ and $D$, the right-hand of (\ref{in1})
gives the correct order of magnitude; examples can be constructed as in the
proofs of Theorems \ref{th2} and \ref{th3} below.

\begin{theorem}
\label{th2}If $G$ is a connected nonbipartite graph of order $n$ and diameter
$D$, then
\begin{equation}
\mu\left(  G\right)  +\mu_{\min}\left(  G\right)  >\frac{2}{\mu^{2D}\left(
G\right)  n}. \label{in2}%
\end{equation}
Moreover, for all $k\geq3,$ $D\geq4,$ and $n=D+2k-1,$ there exists a connected
nonbipartite graph $G$ of order $n$ and diameter $D$ with $\mu\left(
G\right)  >k,$ and
\[
\mu\left(  G\right)  +\mu_{\min}\left(  G\right)  <\frac{4}{\left(
k-1\right)  ^{2D-4}}.
\]

\end{theorem}

Theorem \ref{th2} shows that $\mu\left(  G\right)  +\mu_{\min}\left(
G\right)  $ can be extremely small, although $G$ is nonbipartite and
connected. Here is another viewpoint to this fact.

\begin{theorem}
\label{th3} Let $0<\varepsilon<1/16.$ For all sufficiently large $n,$ there
exists a connected graph $G$ of order $n$ with $\mu\left(  G\right)
+\mu_{\min}\left(  G\right)  <n^{-\varepsilon n}$ such that, to make $G$
bipartite, at least $\left(  1/16-\varepsilon\right)  n^{2}$ edges must be removed.
\end{theorem}

The picture is completely different for regular graphs. In \cite{CGN06} it is
proved that if $G$ is a connected nonregular graph of order $n,$ size $m,$
diameter $D,$ and maximum degree $\Delta,$ then
\[
\Delta-\mu\left(  G\right)  >\frac{n\Delta-2m}{n(D(n\Delta-2m)+1)}.
\]

This result and Theorem \ref{th1} help deduce the following theorems; we omit
their straightforward proofs.

\begin{theorem}
\label{th1.1}If $H$ is a proper subgraph of a connected regular graph $G$ of
order $n$ and diameter $D,$ then%
\[
\mu\left(  G\right)  -\mu\left(  H\right)  >\frac{1}{n(D+1)}.
\]

\end{theorem}

\begin{theorem}
\label{th2.1}If $G$ is a connected regular nonbipartite graph of order $n$ and
diameter $D$, then
\[
\mu\left(  G\right)  +\mu_{\min}\left(  G\right)  >\frac{2}{n(2D+1)}.
\]

\end{theorem}

\begin{theorem}
\label{th4}If $G$ is a connected, nonregular, nonbipartite graph of order $n,$
diameter $D$, and maximum degree $\Delta,$ then%
\[
\Delta+\mu_{\min}\left(  G\right)  >\frac{1}{n(D+1)}+\frac{1}{\mu^{2D}\left(
G\right)  n}.
\]

\end{theorem}

Note that the last two theorems give a fine tuning of a result of Alon and
Sudakov \cite{AlSu00}.

\section{Proofs}

Our proof of Theorem \ref{th1} stems from a result of Schneider \cite{Sch59}
on eigenvectors of irreducible nonnegative matrices; for graphs it reads as:
if $G$ is a connected graph of order $n$ and $x_{\min},$ $x_{\max}$ are
minimal and maximal entries of an eigenvector to $\mu\left(  G\right)  ,$ then%
\[
\frac{x_{\min}}{x_{\max}}\geq\mu^{-n+1}\left(  G\right)  .
\]

We reprove this inequality in a more flexible form that sheds some extra light
on the original matrix result of Schneider as well. Hereafter we write
$dist\left(  u,v\right)  $ for the length of a shortest path joining the
vertices $u$ and $v.$

\begin{proposition}
\label{pro1}If $G$ is a connected graph of order $n$ and $\left(  x_{1}%
,\ldots,x_{n}\right)  $ is an eigenvector to $\mu\left(  G\right)  ,$ then
\begin{equation}
\frac{x_{i}}{x_{j}}\geq\left(  \mu\left(  G\right)  \right)  ^{-dist\left(
i,j\right)  } \label{Sc}%
\end{equation}
for every two vertices $i,j\in V\left(  G\right)  .$
\end{proposition}

\begin{proof}
Clearly we can assume that $i\neq j.$ For convenience we also assume that
$i=1$ and the vertices $\left(  1,\ldots,j\right)  $ form a path joining $1$
to $j.$ Then, for all $u=1,\ldots,j-1,$ we have
\[
\mu x_{u}=\sum_{uv\in E\left(  G\right)  }x_{v}\geq x_{u+1};
\]
hence, (\ref{Sc}) follows by multiplying all these inequalities.
\end{proof}

We shall need also the following simple bound.

\begin{proposition}
\label{pro2}If $G$ is a connected graph of order $n\geq3$ and diameter $D,$
then $\mu^{D}\left(  G\right)  >n/\sqrt{3}.$
\end{proposition}

\begin{proof}
Note that every two vertices can be joined by a walk of $D$ or $D+1$ vertices.
Hence, letting $w_{k}\left(  G\right)  $ be the number of walks of $k$
vertices, we find that $w_{D}\left(  G\right)  +w_{D+1}\left(  G\right)  \geq
n^{2};$ therefore, by a result in \cite{Nik06}, $\mu^{D-1}\left(  G\right)
+\mu^{D}\left(  G\right)  \geq n.$ Since $\mu\left(  G\right)  >\sqrt{2},$ we
see that
\[
\sqrt{3}\mu^{D}\left(  G\right)  >\frac{1}{\sqrt{2}}\mu^{D}\left(  G\right)
+\mu^{D}\left(  G\right)  \geq\mu^{D-1}\left(  G\right)  +\mu^{D}\left(
G\right)  \geq n,
\]
completing the proof.
\end{proof}

\begin{proof}
[\textbf{Proof of Theorem \ref{th1}}]Since $\mu\left(  H\right)  \leq
\mu\left(  H^{\prime}\right)  $ whenever $H\subset H^{\prime}$, we may assume
that $H$ is a maximal proper subgraph of $G$, that is to say, $V\left(
H\right)  =V\left(  G\right)  $ and $H$ differs from $G$ in a single edge
$uv$. Our proof is split into two cases: (a)\emph{ }$H$ connected; (b) $H$
disconnected.\bigskip

\textbf{Case (a): }$H$\emph{ }\textbf{is connected.}

In this case we shall prove a stronger result than required, namely
\begin{equation}
\mu\left(  G\right)  -\mu\left(  H\right)  >\frac{2}{\mu^{2D}\left(  G\right)
n}. \label{in1.1}%
\end{equation}
Our first goal is to prove that, for every $w\in V\left(  H\right)  ,$
\begin{equation}
dist_{H}\left(  w,u\right)  +dist_{H}\left(  w,v\right)  \leq2D. \label{in3}%
\end{equation}
Let $w\in V\left(  H\right)  $ and select in $H$ shortest paths $P\left(
u,w\right)  $ and $P\left(  v,w\right)  $ joining $u$ and $v$ to $w.$ Let
$Q\left(  u,x\right)  $ and $Q\left(  v,x\right)  $ be the longest subpaths of
$P\left(  u,w\right)  $ and $P\left(  v,w\right)  $ having no internal
vertices in common. If $s\in Q\left(  u,x\right)  $ or $s\in Q\left(
v,x\right)  ,$ we obviously have%
\begin{equation}
dist_{H}\left(  w,s\right)  =dist_{H}\left(  w,x\right)  +dist_{H}\left(
s,x\right)  . \label{eq1}%
\end{equation}
The paths $Q\left(  u,x\right)  ,$ $Q\left(  v,x\right)  $ and the edge $uv$
form a cycle in $G;$ write $k$ for its length. Assume that $dist\left(
v,x\right)  \geq dist\left(  u,x\right)  $ and select $y\in Q\left(
v,x\right)  $ with $dist_{H}\left(  x,y\right)  =\left\lfloor k/2\right\rfloor
.$ Let $R\left(  w,y\right)  $ be a shortest path in $G$ joining $w$ to $y;$
clearly the length of $R\left(  w,y\right)  $ is at most $D.$ If $R\left(
w,y\right)  $ does not contain the edge $uv,$ it is a path in $H$ and, using
(\ref{eq1}), we find that%
\begin{align*}
D  &  \geq dist_{G}\left(  w,y\right)  =dist_{H}\left(  w,y\right)
=dist_{H}\left(  w,x\right)  +\left\lfloor k/2\right\rfloor \\
&  =dist_{H}\left(  w,x\right)  +\left\lfloor \frac{dist_{H}\left(
x,u\right)  +dist_{H}\left(  x,v\right)  +1}{2}\right\rfloor \\
&  \geq dist_{H}\left(  w,x\right)  +\frac{dist_{H}\left(  x,u\right)
+dist_{H}\left(  x,v\right)  }{2}=\frac{dist_{H}\left(  w,u\right)
+dist_{H}\left(  w,v\right)  }{2},
\end{align*}
implying (\ref{in3}). Let now $R\left(  w,y\right)  $ contain the edge $uv.$
Assume first that $v$ occurs before $u$ when traversing $R\left(  w,y\right)
$ from $w$ to $y.$ Then
\begin{align*}
dist_{H}\left(  w,u\right)  +dist_{H}\left(  w,v\right)   &  \leq
2dist_{H}\left(  w,x\right)  +dist_{H}\left(  x,u\right)  +dist_{H}\left(
x,v\right) \\
&  \leq2\left(  dist_{H}\left(  w,x\right)  +dist_{H}\left(  x,v\right)
\right)  <dist_{G}\left(  w,y\right)  \leq2D,
\end{align*}
implying (\ref{in3}). Finally, if $u$ occurs before $v$ when traversing
$R\left(  w,y\right)  $ from $w$ to $y,$ then
\begin{align*}
D  &  \geq dist_{G}\left(  w,y\right)  \geq dist_{H}\left(  w,u\right)
+1+dist_{H}\left(  v,y\right) \\
&  =dist_{H}\left(  w,x\right)  +dist_{H}\left(  x,u\right)  +1+dist_{H}%
\left(  v,y\right)  =dist_{H}\left(  w,x\right)  +\left\lceil k/2\right\rceil
\\
&  \geq dist_{H}\left(  w,x\right)  +\frac{dist_{H}\left(  x,u\right)
+dist_{H}\left(  x,v\right)  }{2}=\frac{dist_{H}\left(  w,u\right)
+dist_{H}\left(  w,v\right)  }{2},
\end{align*}
implying (\ref{in3}). Thus, inequality (\ref{in3}) is proved in full.

Let now $\mathbf{x}=\left(  x_{1},...,x_{n}\right)  $ be a unit eigenvector to
$\mu\left(  H\right)  $ and let $x_{w}$ be a maximal entry of $\mathbf{x}$. In
view of (\ref{Sc}) and (\ref{in3}), we have%
\[
\frac{x_{u}x_{v}}{x_{w}^{2}}\geq\frac{1}{\mu^{dist\left(  u,w\right)
+dist\left(  v,w\right)  }\left(  H\right)  }\geq\frac{1}{\mu^{2D}\left(
H\right)  }.
\]
Hence, in view of $x_{w}^{2}\geq1/n,$ we see that
\[
\mu\left(  G\right)  \geq2\sum_{ij\in E\left(  G\right)  }x_{i}x_{j}%
=2x_{u}x_{v}+\mu\left(  H\right)  \geq\frac{2x_{w}^{2}}{\mu^{2D}\left(
H\right)  }+\mu\left(  H\right)  >\frac{2}{\mu^{2D}\left(  G\right)  n}%
+\mu\left(  H\right)  ,
\]
completing the proof of (\ref{in1.1}) and thus of (\ref{in1}).\bigskip

\textbf{Case (b): }$H$\emph{ }\textbf{is disconnected.}

Since $G$ is connected, $H$ is union of two connected graphs $H_{1}$ and
$H_{2}$ such that $v\in H_{1},$ $u\in H_{2}.$ Assume $\mu\left(  H\right)
=\mu\left(  H_{1}\right)  ,$ set $\left\vert H_{1}\right\vert =k,$ and let
$\mathbf{x}=\left(  x_{1},...,x_{k}\right)  $ be a unit eigenvector to
$\mu\left(  H_{1}\right)  .$ Since any maximal entry of $\mathbf{x}$ is at
least $k^{-1/2}$ and $diam$ $H_{1}\leq diam$ $G\leq D,$ Proposition \ref{pro1}
implies that $x_{v}\geq\mu^{-D}\left(  H\right)  k^{-1/2}.$ Set $t=\mu
^{-D}\left(  H\right)  k^{-1/2}$ and consider the unit vector%
\[
\left(  y_{1},...,y_{k},y_{u}\right)  =\left(  x_{1}\sqrt{1-t^{2}}%
,...,x_{k}\sqrt{1-t^{2}},t\right)  .
\]
Then
\begin{align*}
\mu\left(  G\right)   &  \geq\mu\left(  H_{1}+u\right)  \geq2\sum_{ij\in
E\left(  H_{1}+u\right)  }y_{i}y_{j}\geq2t\sum_{uj\in E\left(  H_{1}+u\right)
}y_{j}+2\left(  1-t^{2}\right)  \sum_{ij\in E\left(  H_{1}\right)  }x_{i}%
x_{j}\\
&  \geq2t\sqrt{1-t^{2}}x_{v}+\left(  1-t^{2}\right)  \mu\left(  H\right)
=\frac{1}{\mu^{2D}\left(  H\right)  k}\left(  2\sqrt{1-\frac{1}{\mu
^{2D}\left(  H\right)  k}}-1\right)  +\mu\left(  H\right)  .
\end{align*}
For $k\geq3$, Proposition \ref{pro2} implies that
\begin{align*}
\frac{1}{\mu^{2D}\left(  H\right)  k}\left(  2\sqrt{1-\frac{1}{\mu^{2D}\left(
H\right)  k}}-1\right)   &  >\frac{1}{\mu^{2D}\left(  H\right)  k}\left(
2\sqrt{1-\frac{3}{k^{3}}}-1\right) \\
&  >\frac{1}{\mu^{2D}\left(  H\right)  \left(  k+1\right)  }>\frac{1}{\mu
^{2D}\left(  G\right)  n}.
\end{align*}
Finally, if $k=3,$ then $\mu\left(  H_{1}\right)  =1,$ $\mu\left(  G\right)
\geq\sqrt{2},$ $D\geq2,$ and $n\geq3$; hence,%
\[
\mu\left(  G\right)  -\mu\left(  H\right)  \geq\sqrt{2}-1>\frac{1}{3\left(
\sqrt{2}\right)  ^{4}}\geq\frac{1}{\mu^{2D}\left(  G\right)  n},
\]
completing the proof.
\end{proof}

\bigskip

\begin{proof}
[\textbf{Proof of Theorem \ref{th2}}]Let $\mathbf{x}=\left(  x_{1}%
,...,x_{n}\right)  $ be an eigenvector to $\mu_{\min}\left(  G\right)  $ and
let $V_{1}=\left\{  u:x_{u}<0\right\}  .$ Let $H$ be the maximal bipartite
subgraph of $G,$ containing all edges with exactly one vertex in $V_{1}.$ It
is not hard to see that $H$ is connected proper subgraph of $G$, $V\left(
H\right)  =V\left(  G\right)  ,$ and $\mu_{\min}\left(  H\right)  <\mu_{\min
}\left(  G\right)  .$ Finally, let $H^{\prime}$ be a maximal proper subgraph
of $G$ containing $H.$ We have%
\[
\mu\left(  G\right)  +\mu_{\min}\left(  G\right)  \geq\mu\left(  G\right)
+\mu_{\min}\left(  H\right)  =\mu\left(  G\right)  -\mu\left(  H\right)
\geq\mu\left(  G\right)  -\mu\left(  H^{\prime}\right)  .
\]
and (\ref{in2}) follows from case (a) of the proof of Theorem \ref{th1}.

To construct the required example, set $G_{1}=K_{3}$, $G_{2}=K_{k,k},$ join
$G_{1}$ to $G_{2}$ by a path $P$ of length $n-2k-2$, and write $G$ for the
resulting graph; obviously $G$ is of order $n$ and diameter $n-2k+1$. Set
$\mu=\mu\left(  G\right)  $ and note that $\mu\left(  G\right)  >k.$ Let
$V\left(  G_{1}\right)  =\left\{  u_{1},u_{2},v_{1}\right\}  $ and $P=\left(
v_{1},\ldots,v_{n-2k-1}\right)  ,$ where $v_{n-2k-1}\in V\left(  G_{2}\right)
.$ Let $\mathbf{x}$ be a unit eigenvector to $\mu\left(  G\right)  $ and
assume that the entries $x_{1},x_{2},x_{3},\ldots,x_{n-2k+1}$ correspond to
$u_{1},u_{2},v_{1},\ldots,v_{n-2k-1}.$ Clearly $x_{1}=x_{2},$ and so, from
$\mu x_{2}=x_{2}+x_{3},$ we find that $x_{1}=x_{2}=x_{3}/\left(  \mu-1\right)
.$ Furthermore,%
\[
\mu x_{3}=2x_{2}+x_{4}=\frac{2x_{3}}{\mu-1}+x_{4}<x_{3}+x_{4},
\]
and by induction we obtain $x_{i}<\left(  \mu-1\right)  x_{i+1}$ for all
$3\leq i\leq n-2k.$ Therefore,%
\[
x_{1}=x_{2}\leq\left(  \mu-1\right)  ^{-n+2k+1}x_{n-2k+1}<\left(  k-1\right)
^{-D+2},
\]
and by Rayleigh's principle we deduce that
\[
\mu\left(  G\right)  +\mu_{\min}\left(  G\right)  \leq4x_{1}x_{2}<\frac
{4}{\left(  k-1\right)  ^{2D-4}},
\]
completing the proof.
\end{proof}

\bigskip

\begin{proof}
[\textbf{Proof of Theorem \ref{th3}}]Set $r=\left\lceil n/4\right\rceil +1,$
$s=\left\lceil \left(  1/2-\varepsilon\right)  n\right\rceil ,$ select
$G_{1}=K_{r,r}$, $G_{2}=K_{s},$ join $G_{1}$ to $G_{2}$ by a path $P$ of
length $n-2r-s+1$ and write $G$ for the resulting graph. Note first that, to
make $G$ bipartite, we must remove at least
\[
\binom{s}{2}-\left\lfloor \frac{s^{2}}{4}\right\rfloor \geq\frac{s^{2}}%
{4}-\frac{s}{2}>\frac{\left(  1/2-\varepsilon\right)  ^{2}n^{2}}{4}-\frac
{s}{2}\geq\left(  \frac{1}{16}-\varepsilon\right)  n^{2}%
\]
edges, for $n\ $large enough. Note also that
\[
n-2\left\lceil \frac{n}{4}\right\rceil -2-\left\lceil \left(  \frac{1}%
{2}-\varepsilon\right)  n\right\rceil +1>n-\frac{n}{2}-\left(  \frac{1}%
{2}-\varepsilon\right)  n-4=\varepsilon n-4.
\]
so the length of $P$ is greater than $\varepsilon n-4.$

Let $\mathbf{x}$ be a unit eigenvector to $\mu\left(  G\right)  .$ Clearly the
entries of $\mathbf{x}$ corresponding to vertices from $V\left(  G_{1}\right)
\backslash V\left(  P\right)  $ have the same value $\alpha.$ Like in the
proof of Theorem \ref{th2}, we see that $\alpha<\left(  n/4\right)
^{-\varepsilon n+5}.$ Hence, by Rayleigh's principle, for $n\ $large enough,
we deduce that
\[
\mu\left(  G\right)  +\mu_{\min}\left(  G\right)  \leq4\alpha^{2}\binom{s}%
{2}<\left(  n/4\right)  ^{-2\varepsilon n+10}\frac{n^{2}}{2}<\left(
n/4\right)  ^{-2\varepsilon n+12}<n^{-\varepsilon n},
\]
completing the proof.
\end{proof}

\textbf{Acknowledgment }The author is indebted to B\'{e}la Bollob\'{a}s for
his kind support and to Sebi Cioab\u{a} for interesting discussions.

\end{document}